\documentclass[11pt]{article}

\usepackage{amsmath,amssymb}

\newcommand\zero{{\bf 0}}
\newcommand\x{{\bf x}}
\newcommand\y{{\bf y}}
\newcommand\z{{\bf z}}
\newcommand\bfa{{\bf a}}
\newcommand\bfb{{\bf b}}
\newcommand\bfd{{\bf d}}
\newcommand\bfc{{\bf c}}
\newcommand\bfv{{\bf v}}
\newcommand\bfw{{\bf w}}
\newcommand\bP{{\bf P}}

\newcommand{\bfnu}{{\boldsymbol{\nu}}}

\newtheorem{theorem}{Theorem}[section]
\newtheorem{lemma}[theorem]{Lemma}
\newtheorem{proposition}[theorem]{Proposition}
\newtheorem{definition}[theorem]{Definition}

\newtheorem{corollary}[theorem]{Corollary}

\newcommand{\zz}{{\mathbb Z}}
\newcommand{\qq}{{\mathbb Q}}
\newcommand{\rr}{{\mathbb R}}
\newcommand{\cc}{{\mathbb C}}

\newcommand{\boldf}{{\mathbf f}}
\newcommand{\boldh}{{\mathbf h}}

\newcommand\qed{{\hspace*{\fill}$\Box$\vskip12pt plus 1pt}}

\begin{document}

\title{Polyhedral Methods in Numerical Algebraic Geometry\thanks{Date:
16 October 2008. 
This material is based upon work
supported by the National Science Foundation
under Grant No.\ 0713018.}
}

\author{
Jan Verschelde\thanks{
Department of Mathematics, Statistics, and Computer Science,
University of Illinois at Chicago, 851 South Morgan (M/C 249),
Chicago, IL 60607-7045, USA.
{\em Email:} jan@math.uic.edu or jan.verschelde@na-net.ornl.gov.
{\em URL:} http://www.math.uic.edu/{\~{}}jan.}
}

\date{{\em to Andrew Sommese, on his $60^{th}$ birthday}}

\maketitle

\begin{abstract}
In numerical algebraic geometry witness sets are numerical representations
of positive dimensional solution sets of polynomial systems.
Considering the asymptotics of witness sets we propose certificates
for algebraic curves.  These certificates are the leading terms of
a Puiseux series expansion of the curve starting at infinity.
The vector of powers of the first term in the series is a tropism.  
For proper algebraic curves,
we relate the computation of tropisms to the calculation of mixed volumes.
With this relationship, the computation of tropisms and Puiseux series
expansions could be used as a preprocessing stage prior to a more
expensive witness set computation.  Systems with few monomials have
fewer isolated solutions and fewer data are needed to represent their
positive dimensional solution sets.

\noindent {\bf 2000 Mathematics Subject Classification.}
Primary 65H10.  Secondary 14Q99, 68W30.

\noindent {\bf Key words and phrases.}
certificate,
mixed volume, Newton polytope,
numerical algebraic geometry,
polyhedral method,
polynomial system, 
Puiseux series,
tropism, witness set.

\end{abstract}

\section{Introduction}

Solving polynomial systems numerically used to be restricted to finding
approximations to all isolated solutions.  Via numerical algebraic
geometry~\cite{SW96, SW05}
we are now able to provide numerical representations for all
solutions {\em and} for all dimensions.  
The development of the methods of numerical algebraic geometry 
coincided with upgrades~\cite{SVW03} to the software PHCpack~\cite{Ver99}.
Its blackbox solver computes isolated solutions via polyhedral methods.
The aim of this paper is to bring polyhedral methods to the foreground and
to show how they may lead to finer representations of solution sets.

We consider as
given a polynomial system with at least as many equations as unknowns,
and we want to
determine whether a proper algebraic curve is a solution.
An algebraic curve is proper if it is not contained
in any higher dimensional solution set.  
Moreover, we will assume that the proper algebraic curves we are looking 
for are of multiplicity one.
What we want to determine of a solution curve is a certificate
of its existence.  The certificate should be small and straightforward
to use in computer algebra systems.

For the solution to this problem
we propose polyhedral methods to find
the leading terms of a Puiseux series expansion~\cite{Pui1850}
of the solution curve.
Verifying the existence of a curve as a solution of the given polynomial
system requires then only the formal substitution of the Puiseux series
into the given system, a routine task for computer algebra systems.
The main source of inspiration for this approach was found in the
emerging field of tropical algebraic geometry~\cite{Jen07},
\cite{RST05}, \cite{Stu02}.
To prove that a polynomial system from celestial mechanics
has only finitely many isolated solutions, polyhedral methods 
were applied in~\cite{HM06}.  This paper could be seen as an attempt
to make such proofs automatic.

Working with series expansions is a hybrid form of computation,
combining symbolic and numerical approaches.
Some algebraic curves may actually have finite series expansions
and in the extreme case even consist of only one leading term.
When considering systems with natural parameters, 
such solution curves may be more useful than other solutions
and it may be worthwhile to look for such solutions first.
Our goal is to develop a polyhedral method that will focus
on computing one dimensional solutions.

After defining witness sets and outlining our problem statement,
we will introduce tropisms for binomial systems.
Proper algebraic curves defined by binomial systems admit
a very explicit solution.  To examine the normalization of
tropisms we consider deformations of witness sets.
The link between tropisms and mixed volumes passes through
the second theorem of Bernshte\v{\i}n~\cite{Ber75}.
We show how to compute
tropisms to proper algebraic curves defined by square systems
via a special lifting.

In the spirit of the theme of the conference, held in honor of Andrew
Sommese, 22-24 May 2008, this paper is on the interactions of 
classical algebraic geometry --- in particular the Puiseux series,
the roots of polyhedral methods and tropical algebraic geometry ---
and the new field of numerical algebraic geometry.
Detailed descriptions of algorithms and their implementations
are still a work in progress.

\noindent {\bf Acknowledgements.}  Preliminary versions of the ideas
in this paper were presented at four conferences held respectively
in Tobago, Notre Dame, Hong Kong, and Vancouver.  
The author thanks all organizers of these meetings for their invitations
and the participants for their feedback.

\section{Witness Sets and Newton Polytopes}

In numerical algebraic geometry, 
homotopy continuation methods~(\cite{Li03}, \cite{Mor87}, \cite{SW05})
manipulate algebraic sets.  The key data representation for
algebraic sets is a witness set (first defined as witness point
set in~\cite{SVW01}, see also~\cite{SW05}), defined below.

\begin{definition} \label{defwitset} {\rm
Given a system $\boldf(\x) = \zero$, we represent
a component of $\boldf^{-1}(\zero)$ of dimension $k$ and degree $d$ by
{\em a witness set} which consists of $\boldf$,
$k$ general hyperplanes $L$ and
$d$ generic points in $\boldf^{-1}(\zero) \cap L$. }
\end{definition}
A generic choice of coefficients for the $k$ hyperplanes implies
that all solutions will be isolated and regular, 
unless the component has a multiplicity higher than one.
Definition~\ref{defwitset} incorporates the theorem of Bertini,
see~\cite{SW05} for its application in numerical algebraic geometry.
Bertini~\cite{BHSW08} is a more recent software system
for numerical algebraic geometry.
A dictionary compares in~\cite{SVW08} witness sets
to lifting fibers~\cite{GH01}, \cite{GLS01} 
in a geometric resolution~\cite{GHMP95}.

We call a system square if it has as many equations as unknowns.
After adding linear equations to a square system, we~\cite{SV00} 
add slack variables in an embedding to make the overdetermined system 
square again.
Using a flag of linear spaces, defined by an decreasing sequence
of subsets of the $k$ general hyperplanes,
\begin{equation}
   L = L_k \supset L_{k-1} \supset \cdots \supset L_1 \supset L_0 = \emptyset,
\end{equation}
we move solutions with nonzero slack values to
generic points on lower dimensional components,
using a cascade of homotopies~\cite{SV00}.
By default, the top dimension~$k$ in the cascade starts at~$n-1$,
but if we know the system has only solution curves we start with $k=1$.
According to~\cite{GH93}, the dimension of an algebraic variety
can be determined in polynomial time.

The cost of the homotopy cascade algorithm is determined by the number
of paths, starting at the solutions of the top dimensional system.
By the embedding, this top dimensional system has only isolated roots
and thus one may apply any solver to compute those isolated roots.
To count the number of solutions for square systems with generic
coefficients we look at the Newton polytopes.
The Newton polytopes are spanned by the exponents of the monomials
which occur with nonzero coefficients in the system.
We formally describe sparse polynomials as follows.

\begin{definition}  {\rm Given a polynomial
\begin{equation}
   f(\x) = \sum_{\bfa \in A} c_\bfa \x^\bfa, \quad c_\bfa \not= 0,
   \quad \x^\bfa = x_1^{a_1} x_2^{a_2} \cdots x_n^{a_n},
\end{equation}
the {\em support} $A$ of~$f$ collects only those exponent vectors
with nonzero coefficient in~$f$.  
The convex hull of~$A$ is the {\em Newton polytope} of~$f$.
We denote the inner product by $\langle \cdot , \cdot \rangle$.
To define faces ${\partial}_\bfv P$ of~$P$ 
we use a {\em support function} $p$:
\begin{equation} \label{eqsupfun}
    p(\bfv) = \min_{\x \in P} {\langle \x , \bfv \rangle}
    \quad {\rm so} \quad
    {\partial}_\bfv P
     = \{ \ \x \in P \ | \ \langle \x , \bfv \rangle = p(\bfv) \ \}.
\end{equation}
The equation $\langle \x , \bfv \rangle = p(\bfv)$ determines
a {\em supporting hyperplane} for the face~$\partial_\bfv P$.
A vector~$\bfv$ perpendicular to
a $k$-dimensional face ${\partial}_\bfv P$ of $P$
lies in an $(n-k)$-dimensional cone. }
\end{definition}

In defining support functions as in~(\ref{eqsupfun}),
we choose the minimum instead of
the maximum convention, opting for inner rather than outer normals.
These minimum and maximum conventions correspond to letting a 
deformation parameter respectively go to zero or to infinity.
In~\cite{GKZ94}, Newton polytopes arise as compactifications of
amoebas~\cite{Mik04}, obtained by taking logarithms of 
the variety~\cite{Ber71}.
Algorithms to compute amoebas are presented in~\cite{The02}.
We refer to~\cite{Tho06} and~\cite{Zie95} for references on polytopes.
For a system $\boldf(\x) = \zero$,
we collect the Newton polytopes of the tuple of polynomials
$\boldf = (f_1,f_2,\ldots,f_n)$ 
in the tuple $\bP = (P_1,P_2,\ldots,P_n)$.

\begin{definition} {\rm Given a tuple of Newton polytopes
${\bP} = (P_1,P_2,\ldots,P_n)$ we define {\em the mixed volume}
$V_n({\bP})$ via the formula
\begin{equation} \label{eqmixvol}
  V_n (P_1,P_2,\ldots,P_n) =
  \sum_{\begin{array}{c}
            \bfv \in \zz^n \\ \gcd(\bfv) = 1
        \end{array} } \ p_1 ({\bf v}) \
  V_{n-1}({\partial}_\bfv P_2, \ldots , {\partial}_\bfv P_n),
\end{equation}
where $p_1$ is the support function for~$P_1$.
Vectors $\bfv$ are normalized so the components of~$\bfv$
have their greatest common divisor equal to 1. }
\end{definition}
Because polytopes are spanned by only finitely many points,
only finitely many vectors $\bfv$ will yield a nonzero
contribution to~(\ref{eqmixvol}).
In the special case: $P = P_1 = P_2 = \cdots = P_n$: 
$V_n({\bP}) = n! {\rm volume}(P)$.
Handling the mixed case, when several but not all polytopes are
repeated, is important for efficient algorithms,
developed in~\cite{GL03}, \cite{LL08}, and~\cite{MTK07}.
Mixed volumes are at the core of classical geometry~\cite{Sch93}
and have applications to tomography~\cite{Gar06}.

\begin{theorem} [Bernshte\v{\i}n Theorem A~\cite{Ber75}] \label{theobera}
The number of roots of a generic system equals the mixed volume of 
its Newton polytopes.  For any system, the mixed volume bounds the
number of isolated solutions in~$(\cc^*)^n$,
$\cc^* = \cc \setminus \{ 0 \}$.
\end{theorem}

In the same paper, Bernshte\v{\i}n gave in his second theorem
precise conditions for the mixed volume to be sharp.
Because of these conditions, we speak of {\em count}~\cite{CR91} 
instead of {\em bound}, see also~\cite{Roj99, Roj03}.
We will formulate the second theorem later,
but now we have introduced enough terminology to formulate
our problem statement.
The proof for Theorem~\ref{theobera} given in~\cite{Ber75}
is constructive and served as basis for polyhedral method
in~\cite{VVC94} to find all isolated solutions.  A more
general polyhedral homotopy homotopy method was developed 
in~\cite{HS95}, see also~\cite{Stu98}.
The authors of~\cite{JMSW08} address the 
complexity of Bernshte\v{\i}n's first theorem.
The development of recent software for polyhedral homotopies is
described in~\cite{GKKTFM04} and~\cite{LLT08}.

We list three problems with the current use of witness sets.
First and foremost,
as we add hyperplanes to the system, the polyhedral root count
increases for sparse systems.
For example, for a benchmark problem like the cyclic 8-roots system,
adding one hyperplane raises the mixed volume from 2,560 to 4,176.
For cyclic 12-roots, the same operation brings the mixed volume
from 500,352 to 983,952.  While the cascade will give us all
solutions at the end, if we are only interested in the curves,
we do not want to compute start solutions to the isolated roots.
The second problem concerns symmetry.  Many polynomial systems
have obvious permutation symmetries and with
polyhedral methods we can setup symmetric homotopies~\cite{VG95}
to compute only the generators of each orbit of isolated roots.
However, extending the symmetric polyhedral homotopies of~\cite{VG95}
to deal with positive dimensional solution sets conflicts
with the current witness set data representation.
Thirdly, users of numerical algebraic geometry methods require
and need guarantees for the results to be correct.
Once we have a witness set, any path tracker may be used to sample
points on the solution set.  While path trackers are standard in
numerical analysis, users unfamiliar or uncomfortable with 
floating-point computations require exact answers.
Additionally, to determine the degree of a solution set correctly,
we must be able to certify that all solutions have been found.

\section{Proper Algebraic Curves defined by Binomial Systems}

A binomial system has exactly two monomials
with a nonzero coefficient in every equation.
We consider $n-1$ equations in $n$ variables.
By limiting the number of equations and restricting to
nontrivial solutions in~$(\cc^*)^n$ we are reducing the complexity
of the problem.
However, even already for binomial systems,
the complexity of counting all isolated solutions
is \#$P$-complete~\cite{CD07}.

For example ($n = 3$):
\begin{equation} \label{eqbinsys}
   \left\{
      \begin{array}{r}
          x_1 x_2^2 x_3 - 2 x_1^2 x_2^3 x_3 = 0 \\ \vspace{-4mm} \\
        3 x_1^2 x_2^2 x_3^5 + 9 x_1 x_2 x_3 = 0 \\
      \end{array}
   \right.
   \quad \equiv \quad
   \left\{
      \begin{array}{lllcr}
          x_1^{-1} \! & \! x_2^{-1} \! & \! \!       & \! = \! & \! 2
 \\ \vspace{-4mm} \\
          x_1^2    \! & \! x_2      \! & \! x_3^4 \! & \! = \! & \! -3
      \end{array}
   \right.
\end{equation}
The system at the right of~(\ref{eqbinsys}) is a normal form
of the system.  For general $n$, we can always write a binomial system 
as a tuple of equations of the form $\x^\bfa = c$.
Writing a system in this normal form removes trivial solutions 
with zero coordinates.
For a binomial system in its normal form, we collect all exponents
in a matrix $A \in \zz^{(n-1) \times n}$ and 
its coefficients in $\bfc \in (\cc^*)^{n-1}$.
Continuing the example, we have:
\begin{equation}
   A = 
  \left[
    \begin{array}{rrr}
       -1 & -1 & 0 \\
        2 & 1 & 4
    \end{array}
  \right],
  \quad {\rm rank}(A) = 2,
  \quad
  \bfv =
  \left[
     \begin{array}{r}
        4 \\ -4 \\ -1
     \end{array}
  \right]:
  A \bfv = \zero.
\end{equation}
The vector~$\bfv$ in the kernel of~$A$ will determine the shape
of the solution curve.

For $A$: ${\rm rank}(A) = n-1$, there is a unique vector~$\bfv$ 
in the kernel.  
A binomial system of $n-1$ equations in~$n$ variables will have
have a proper (i.e.: not contained in any other
higher dimensional solution set) solution curve in~$\cc^*$ if
and only if the rank of the exponent matrix $A$ is~$n-1$. 
Although a normal form of a binomial system is
not unique, the vector~$\bfv$ for a proper algebraic curve is
independent of the choice of a particular normal form.
We will show that proper algebraic curves have solutions of the 
type $x_k = c_k t^{v_k}$, with $c_k \in \cc^*$, $k=1,2,\ldots,n$.

To simplify the system,
we use~$\bfv$ to define a unimodular matrix~$M$ ($\det(M) = 1$) and
a coordinate transformation (called a power transformation in~\cite{Bru00})
denoted by $\x = \y^M$:
\begin{equation}
   M =
   \left[
      \begin{array}{rrr}
         +4 & 0 & 1 \\
         -4 & 1 & 0 \\
         -1 & 0 & 0
      \end{array}
   \right]
   \quad
   AM = 
   \left[
      \begin{array}{rrr}
         0 & -1 & -1 \\
         0 & 1 & 2
      \end{array}
   \right]
   \quad
   \left\{
      \begin{array}{l}
         x_1 = y_1^{+4} y_3 \\ \vspace{-4mm} \\
         x_2 = y_1^{-4} y_2 \\ \vspace{-4mm} \\
         x_3 = y_1^{-1}
      \end{array}
   \right.
\end{equation}
After applying the coordinate transformation, defined by~$M$,
we can divide out the variable~$y_1$ and we obtain
a system of two equations in two unknowns: $y_2$ and $y_3$.
For our example, we find one solution
\begin{equation}
   \begin{array}{c}
      (y_2 = -\frac{1}{12}, y_3 = -6) 
   \end{array}
   \quad {\rm to} \quad
   \left\{
      \begin{array}{llcr}
         y_2^{-1} \! & \! y_3^{-1} \! & \! = \! &  2~\! 
\\ \vspace{-4mm} \\
         y_2      \! & \! y_3^2    \! & \! = \! & -3.
      \end{array}
   \right.
\end{equation}
We rename the free variable~$y_1$ to~$t$ and use the values
found for~$y_2$ and $y_3$ in the representation of the solution curve,
see the left of~(\ref{eqdegsol}).
To compute the degree of the curve,
we add a random hyperplane and substitute the expression
found for the solution:

\begin{equation} \label{eqdegsol}
   \left\{
      \begin{array}{l}
         x_1 = -6 t^{+4} \\ \vspace{-4mm} \\
         x_2 = -\frac{1}{12} t^{-4} \\ \vspace{-4mm} \\
         x_3 = t^{-1}
      \end{array}
   \right.
   \quad \quad
   \left\{
      \begin{array}{lllcrc}
          x_1^{-1} \! & \! x_2^{-1} \! & \! \!       & \! = \! & \! 2 &
 \\ \vspace{-4mm} \\
          x_1^2    \! & \! x_2      \! & \! x_3^4 \! & \! = \! & \! -3 &
 \\ \vspace{-4mm} \\
          \multicolumn{6}{c}{\gamma_0 + \gamma_1 x_1
              + \gamma_2 x_2 + \gamma_3 x_3 = 0}
      \end{array}
   \right.
\end{equation}
The coefficients $\gamma_0$, $\gamma_1$, $\gamma_2$, and $\gamma_3$
are random complex numbers.
After substitution the left of~(\ref{eqdegsol})
into the right of~(\ref{eqdegsol}) we obtain
$\gamma_0 + \gamma_1 ( -6 t^{+4} ) 
 + \gamma_2 ( -\frac{1}{12} t^{-4} ) + \gamma_3 t^{-1}$.
Clearing denominators, we find a polynomial in~$t$ of degree~8.

Given the leading term of the Puiseux series expansion for an
algebraic curve, the degree of the curve follows from
the leading exponents of the Puiseux series.
We can formalize this in the following proposition:
\begin{proposition} \label{propdegbinsys}
Consider a proper algebraic curve defined by $\x^A = \bfc$,
with $\bfv$: $A \bfv = \zero$.  Let $M$ be the unimodular matrix
in the transformation $\x = \y^M$ that eliminates~$y_1$.
Let $B$ be the matrix obtained from removing the first zero column
of~$AM$.
Then the degree of the curve equals
\begin{equation} \label{eqdegformula}
   |\det(B)| \times | \max_{i=1}^n v_i - \min_{i=1}^n v_i |.
\end{equation}
\end{proposition}

\noindent {\em Proof.}  The form of the solution curve is
$x_k = \alpha_k t^{v_k}$ where the coefficients $\alpha_k$ are
the nonzero roots of $\y^{AM} = \bfc$.
The number of roots equals $|\det(B)|$.
The curve will have as many components as~$|\det(B)|$.
Because every component has the same degrees in~$t$,
every component will have the same degree, so it suffices
to compute the degree of one component.

To compute the degree of one component, we reduce it to as many
isolated points as its degree, intersecting it with 
a hyperplane with random coefficients.
Substituting the form of the solution component
into that hyperplane yields a polynomial in~$t$.
To clear denominators we multiply by the most negative
exponent or we divide out trivial solutions by subtracting the
lowest positive exponent of~$\bfv$.~\qed

We point out that the formula~(\ref{eqdegformula}) is {\em not}
invariant to unimodular transformations.
Consider for example the plane curve defined by 
$f(x_1,x_2) = x_1 x_2 - 1 = 0$ with $\bfv = (+1,-1)$.
The unimodular coordinate transformation defined by~$\bfv$: 
$x_1 = y_1$, $x_2 = y_1^{-1} y_2$ reduces $f$ to $y_2 - 1 = 0$.
In the new $y$-coordinates, we now have a line $(y_1,1)$.

For a proper algebraic curve defined by a binomial system
we can provide an exact certificate for its degree,
{\em independently} of the choice of the coefficients.
The vector~$\bfv$ we computed in the example above is an example
of a tropism.  We define tropisms for general polynomial systems
in the next section.

\section{Tropisms and Initial Forms}

Puiseux series occupy a central role in the study of algebraic curves.
The leading exponents of the series are called tropisms, defined by
the Newton polytopes of a polynomial system.

\begin{definition} {\rm
Consider $\boldf(\x) = \zero$ with
Newton polytopes in $(P_1,P_2,\ldots,P_N)$.
A {\em tropism} is a vector perpendicular
to one edge of each $P_i$, for $i=1,2,\ldots,N$. }
\end{definition}

Our definition of tropisms differs from the usual one in the literature
(\cite{LTR08}, \cite{Mau80}) where all coordinates in the tropism
(or critical tropism in singularity theory~\cite{LTR08})
are required to be positive.
This requirement is natural if one looks for one point on a
solution curve in affine space and then shifts that point to the origin.
The only normalization we will require is on the sign of the first
coordinate.  We will provide arguments for this normalization when
we consider the asymptotics of witness sets in the next section.
The other difference in the definition is then that we may have tropisms
pointing out isolated solutions at infinity, solutions that do not
give rise to an initial term in a Puiseux series expansion.
Using a more refined terminology we could give the $\bfv$
in Definition~\ref{deftropism} the name {\em pretropism},
as a part of a tropical prevariety.
A tropical prevariety~\cite{BJSST07} corresponds to the intersection of 
the normal cones of the polytopes.

The edges perpendicular to a tropism are Newton polytopes 
of {\em an initial form system} which may have solutions 
in $(\cc^*)^n$.

\begin{definition} \label{deftropism} {\rm
Let $\bfv \in \zz^n \setminus \{ 0 \}$
and ${\displaystyle f(\x) = \sum_{\bfa \in A} c_\bfa \x^\bfa}$.
Denoting the inner product by $\langle \cdot , \cdot \rangle$,
the {\em initial form of $f$ in the direction $\bfv$} is
\begin{equation}
{\rm in}_\bfv f(\x)
= \sum_{\begin{array}{c}
           \bfa \in A \\ \langle \bfa, \bfv \rangle = m
        \end{array}} c_\bfa \x^\bfa
\quad {\rm with} \quad
m = \min \{ \ \langle \bfa, \bfv \rangle \ | \ \bfa \in A \ \}.
\end{equation}
Let $\boldf = (f_1,f_2,\ldots,f_N)$ be a tuple of polynomials.
For $\bfv \in \zz^n \setminus \{ 0 \}$, the {\em initial form system}
${\rm in}_\bfv \boldf(\x) = \zero$ is defined by the tuple 
${\rm in}_\bfv \boldf = ( {\rm in}_\bfv f_1, {\rm in}_\bfv f_2, 
\ldots,  {\rm in}_\bfv f_N).$
}
\end{definition}

Although we define initial form systems via tropisms,
vectors perpendicular to edges of the Newton polytopes,
the number of points in the support of each initial form
is at least --- but not exactly --- two.
Although the initial form system is thus not necessarily
a binomial system, we can always eliminate one variable.

In~\cite{Bru00} and ~\cite{Kaz99},
systems supported on faces of Newton polytopes are
called truncated systems.  This terminology refers to the process of
substitution the power series and then selecting those terms that
correspond to the lowest power in the variable of the series.
We prefer to use initial forms because of the relationship with
term orders widely used for Gr\"obner bases~\cite{Stu96}, \cite{Tho06}.

Of special importance are the conditions for which
the mixed volume is sharp, formulated in the tropical language.

\begin{theorem} [Bernshte\v{\i}n Theorem B~\cite{Ber75}]
Consider $\boldf(\x) = \zero$, $\boldf = (f_1,f_2,\ldots,f_n)$, 
$\x = (x_1,x_2,\ldots,x_n)$.  
If for all tropisms $\bfv$: ${\rm in}_\bfv \boldf(\x) = \zero$ has
no solutions in $(\cc^*)^n$, 
then $\boldf(\x) = \zero$ has exactly
as many isolated solutions in $(\cc^*)^n$ as $V_n({\bP})$.
\end{theorem}

If there is a positive dimensional solution set, then this set
stretches out to infinity and the system $f(\x) = \zero$ must
have solutions at infinity.  
Of particular interest is the following consequence
of Bernshte\v{\i}n's second theorem.

\begin{corollary}
If there are no tropisms $\bfv$
for which ${\rm in}_\bfv \boldf(\x) = \zero$ has roots in $(\cc^*)^n$
then $\boldf(\x) = \zero$ has no solutions at infinity.
\end{corollary}

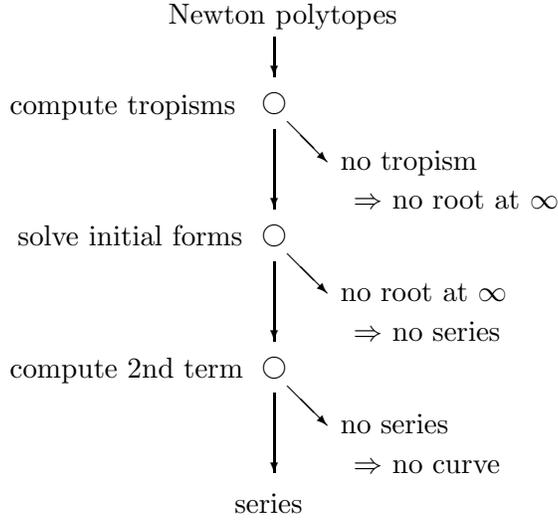
\begin{figure}[ht]
\begin{center}
\begin{picture}(200,180)(0,20)
\put(60,190){Newton polytopes}
\put(0,156){compute tropisms}
\put(100,185){\vector(0,-1){15}} \put(100,160){\circle{8}}
\put(105,153){\vector(+1,-1){15}}
\put(125,135){no tropism}
\put(130,120){$\Rightarrow$ no root at $\infty$}

\put(3,106){solve initial forms}
\put(100,150){\vector(0,-1){30}} \put(100,110){\circle{8}}
\put(105,103){\vector(+1,-1){15}}
\put(125,85){no root at $\infty$}
\put(130,70){$\Rightarrow$ no series}

\put(0,56){compute 2nd term}
\put(100,100){\vector(0,-1){30}}  \put(100,60){\circle{8}}
\put(105,53){\vector(+1,-1){15}}
\put(125,35){no series}
\put(130,20){$\Rightarrow$ no curve}
\put(100,50){\vector(0,-1){30}}
\put(85,5){series}
\end{picture}
\end{center}
\label{figstages}
\caption{Computing a certificate for a proper algebraic curve.}
\end{figure}

The tropisms and the roots of the
corresponding initial form systems give the leading term of the
Puiseux series expansion of the solution curve.
The second term in the Puiseux series expansion will provide
the certificate for the existence of the solution curve.
Also the coefficients of the second term will be roots of
a truncated polynomial system.  Once we have those roots,
we can further grow the Puiseux series expansion symbolically,
or apply numerical predictor-corrector methods to sample points
along the solution curve.
In Figure~\ref{figstages} (slightly adapted from~\cite{AV08})
we sketch the idea for computing this certificate.

The computation of the second term in the Puiseux series
goes along the so-called Newton-Puiseux method, as outlined
in the proof of the theorem of Puiseux~\cite{Wal50},
see also~\cite{dJP00}.
In~\cite{AMNR92} algorithms for Puiseux series for space curves
are described and an implementation in CoCoA is mentioned.
General fractional power series solutions are described in~\cite{McD02}.
See~\cite{JMM08}, \cite{JMSW08} and~\cite{PR08}
for recent symbolic algorithms, and~\cite{Pot07}, \cite{PR08b} for
a symbolic-numeric approach.

Algebraic curves defined by binomial systems can be solved by one tropism
but it may happen that one tropism solves a more general system.
Consider for example the cyclic 4-roots problem.

\begin{equation} \label{eqcyclic4}
  \boldf(\x) = \left\{
  \begin{array}{c}
     x_1 + x_2 + x_3 + x_4 = 0 \\
     x_1 x_2 + x_2 x_3 + x_3 x_4 + x_4 x_1 = 0 \\
     x_1 x_2 x_3 + x_2 x_3 x_4 + x_3 x_4 x_1 + x_4 x_1 x_2 = 0 \\
     x_1 x_2 x_3 x_4 - 1 = 0 \\
  \end{array}
\right.
\end{equation}
There is one tropism $\bfv = (+1,-1,+1,-1)$ which leads
to the initial form system ${\rm in}_\bfv f(\z) = \zero$:

\begin{equation}
{\rm in}_\bfv \boldf(\x) = \left\{
  \begin{array}{c}
     x_2 + x_4 = 0 \\
     x_1 x_2 + x_2 x_3 + x_3 x_4 + x_4 x_1 = 0 \\
     x_2 x_3 x_4 + x_4 x_1 x_2 = 0 \\
     x_1 x_2 x_3 x_4 - 1 = 0 \\
  \end{array}
\right.
\quad
\left\{
   \begin{array}{l}
      x_1 = y_1^{+1} \\
      x_2 = y_1^{-1} y_2 \\
      x_3 = y_1^{+1} y_3 \\
      x_4 = y_1^{-1} y_4 \\
   \end{array}
\right.
\end{equation}
The system ${\rm in}_\bfv \boldf(\y) = \zero$ has two solutions.
and we find the two solution curves:
$\left( t, -t^{-1}, -t, t^{-1} \right)$ and
$\left( t, t^{-1}, -t, -t^{-1} \right)$.

Note that $\bfv = (-1,+1,-1,+1)$ is a tropism as well
for the solution curves of cyclic 4-roots, but considering
this tropism corresponds to setting $x_1 = t^{-1}$ or moving
the curve to infinity instead of to zero as $t$ goes to zero.
We will examine this in greater generality in the next section.

\section{Asymptotics of Witness Sets}

One way to compute tropisms would be to start from a witness
set for an algebraic curve in $n$-space given by $d$ points on a
general hyperplane $c_0 + c_1 x_1 + c_2 x_2 + \cdots + c_n x_n = 0$
and satisfying a system $\boldf(\x) = \zero$.
We then deform a witness set for a curve in two stages:
\begin{enumerate}
\item The first homotopy moves to a hyperplane in special position:
\begin{equation} \label{eqasyhom1}
   \boldh(\x,t) = 
   \left\{
      \begin{array}{cl}
          \boldf(\x) = \zero \\
          (c_0 + c_1 x_1 + \cdots + c_n x_n) t
          + (c_0 + c_1 x_1) (1-t) = 0,
          & {\rm for~} t {\rm ~from~1~to~0.}
      \end{array}
   \right.  
\end{equation}
\item After renaming $c_0 + c_1 x_1 = 0$ into $x_1 = \gamma$,
      we let $x_1$ go to zero with the following homotopy:
\begin{equation} \label{eqasyhom2}
   \boldh(\x,t) =
   \left\{
      \begin{array}{cl}
          \boldf(\x) = \zero \\
          x_1 - \gamma t = 0,
          & {\rm for~} t {\rm ~from~1~to~0.}
      \end{array}
   \right. 
\end{equation}
\end{enumerate}

The two homotopies need further study.  
In the first homotopy~(\ref{eqasyhom1}) some paths will diverge, 
consider for example $f(x_1,x_2) = x_1 x_2 - 1$.
Even all paths may diverge if the solution curve lies in some
hyperplane perpendicular to the first coordinate axis $x_1 = c$
with $c$ different from $-c_0/c_1$.

\begin{lemma} \label{lemasyhom1}
All solutions at the end of the homotopy $\boldh(\x,t) = \zero$
of~{\rm (\ref{eqasyhom1})} lie on the curve defined by $\boldf(\x) = \zero$
and in the hyperplane~$x_1 = -c_0/c_1$.
\end{lemma}

\noindent {\em Proof.}  We claim that we find the same solutions
to $\boldh(\x,t=0) = \zero$ either by using the homotopy in~(\ref{eqasyhom1})
or by solving $\boldh(\x,0) = \zero$ directly.
This claim follows from cheater's homotopy~\cite{LSY89} 
or the more general
coefficient-parameter polynomial continuation~\cite{MS89}.~\qed

The second claim we make in the main theorem below is that
we recover all data lost with tropisms.
For simplest example of the hyperbola $x_1 x_2 - 1 = 0$:
its solution is~$(x_1 = t, x_2 = t^{-1})$
and the tropism is $\bfv = (1,-1)$.
The lemma below extends the normal form for the Puiseux
series expansion for plane curves (as used in~\cite{Wal50})
to general space curves.

\begin{lemma} \label{lemasyhom2}
As $t \rightarrow 0$ in the homotopy~{\rm (\ref{eqasyhom2})},
the leading powers of the Puiseux series expansions are the
components of a tropism.  In particular, the expansions have the form
\begin{equation} \label{eqformsol}
   \left\{
     \begin{array}{l}
        x_1 = t \\
        x_k = c_k t^{v_k}(1+O(t)), \quad k=2,\ldots,n.
     \end{array}
   \right.
\end{equation}
\end{lemma}

\noindent {\em Proof.} Following Bernshte\v{\i}n's second theorem,
a solution at infinity is a solution in~$(\cc^*)^n$ of an initial
form system.  For a solution to have values in~$(\cc^*)^n$, 
all equations in that system need to have at least two monomials.
So the system is an initial form system defined by a tropism.
To arrive at the form of~(\ref{eqformsol}) for the solution
defined by the homotopy~(\ref{eqasyhom2}) we rescale the
parameter~$t$ so we may replace $x_1 = \gamma t$ by $x_1 = t$.~\qed

Also in the second homotopy, solution paths are most likely to
diverge and the directions of the diverging paths are defined
by the tropisms.  The numerical computation of those directions
can be done by endgames using power series as 
in~\cite{MSW92} and~\cite{HV98}.

\begin{definition} \label{defcertificate} {\rm
Given a system $\boldf(\x) = \zero$ which defines a proper algebraic curve.
Consider a Puiseux series expansion of the form
\begin{equation} \label{eqcertificate}
  \left\{
     \begin{array}{rcl}
        x_1 & = & t \\
        x_k & = & c_k t^{v_k} + d_k t^{w_k} + \cdots
          \quad c_k, d_k \in \cc^*, 
          v_k, w_k \in \qq, k=2,\ldots,n.
     \end{array}
  \right.
\end{equation}
Then {\em a certificate} for the solution curve consist of
the exponents $(v_2,\ldots,v_n)$, $(w_2,\ldots,w_n) \in \qq^{n-1}$ and
the coefficients $(c_2,\ldots,c_n)$, $(d_2,\ldots,d_n) \in (\cc^*)^{n-1}$.
}
\end{definition}

The tropism~$\bfv$ (or pretropism) shows there are solutions at infinity,
but solutions at infinity could be isolated.  
In that case the tropism~$\bfv$ is still a certificate (but then more 
like a death certificate) for the lack of sharpness of the mixed volume
to count all isolated roots.
The exponents~$\bfw$ of the second term in the series show the solution 
is part of a curve.
Arguing in favor of extending the term tropism rather than reserving
it only for the leading exponents in the series, we point out that
a certificate consists of a pair of vectors: ($\bfv$, $\bfw$)
and a corresponding pair of solutions $(\bfc,\bfd)$, respectively
of an initial form system and a truncated system.
As solutions of overdetermined polynomial systems, 
the coefficients $(\bfc,\bfd)$ can be certified~\cite{DS99}
by $\alpha$-theory~\cite{BCSS98}.

Relating the data ($\bfv$, $\bfw$) and $(\bfc,\bfd)$ to
a witness set, we note that the data corresponds geometrically
to cutting the curve with a special hyperplane $x_1 = 0$.
There is the risk of missing curves in coordinate hyperplane
$x_1 = c$, for $c \in \cc^*$ and there may also be singular
solutions occurring for $x_1 = 0$.  Except for these two difficulties, 
the certificate provides a predictor to sample the solution curve.
The capability to sample a solution curve is intrinsic in the
definition of a witness set, but to verify this capability one
needs to compute at least one step of Newton's method at at least
one point of the witness set.  Substituting the certificate in
the original system --- with and without the second term and watching
the degree in~$t$ increase --- is a more elementary operation.

On the one hand,
we can view the certificate as a very special witness set,
obtained by intersecting the curve with the hyperplane~$x_1 = 0$.
On the other hand,
we can view the certificate as a very special lifting fiber,
where one free variable $x_1$ is specialized to zero.

\begin{lemma}  The certificate of~{\rm (\ref{eqcertificate})} may
be written in the form
\begin{equation}
   \left\{
      \begin{array}{lcll}
         x_1 & = & t^{\nu_1}, & \nu_1 \geq 1, \\
         x_k & = & \alpha_k t^{\nu_k} + \beta_k t^{\mu_k}, & k=2,\ldots,n,
      \end{array}
   \right.
\end{equation}
where $\alpha_k, \beta_k \in \cc^*$ and $\nu_k, \mu_k \in \zz$
and $\nu_1$ is the smallest natural number to clear the
denominators in the series for the other components $x_k$, $k > 1$.
Then the degree of the branch is determined by
\begin{equation}
   \#R \times | \max_{i=1}^n \nu_i - \min_{i=1}^n \nu_i | 
\end{equation}
where $R$ is the set of initial roots of the initial form system
${\rm in}_\bfnu \boldf(\x = \y^M) = \zero$,
for $M$ a unimodular matrix with first 
column equal to~$\bfnu = (\nu_1,\nu_2,\ldots,\nu_n)$.
\end{lemma}

\noindent {\em Proof.}  The existence of~$\nu_1$ follows from
the definition of Puiseux expansions.  Suppose there would not be
a smallest~$\nu_1$ to clear denominators, then we could make
a plane curve for which we could not clear denominators.
Proposition~\ref{propdegbinsys} is generalized in two ways:
\begin{enumerate}
\item The initial form system ${\rm in}_\bfnu \boldf(\y) = \zero$
      is no longer a binomial system for which we can count
      the number of roots via a determinant.
      Instead we now solve ${\rm in}_\bfnu \boldf(\y) = \zero$ and
      collect the roots in~$R$.
\item We have to show that only the leading terms in the Puiseux
      series expansion determine the degree of the solution branch.
      Consider $t \rightarrow 0$ and consider the solutions
      in the hyperplane $x_1 = t^{\nu_1}$.  By the form of the expansion
      we cannot have more solutions for $t > 0$ then we would have 
      at $t = 0$.  At $t = 0$, only the leading terms matter in
      defining the initial form 
      system~${\rm in}_\bfnu \boldf(\y) = \zero$.~\qed
\end{enumerate}

The key point of the two homotopies we considered above
is the argument for the normalization
of the first coordinate $v_1$ of the tropisms to $v_1 > 0$.
Moving the degree many points on the curve on
the hyperplane $x_1 = \gamma$ to $x_1 = 0$ with $x = t$
is equivalent to moving those points to infinity with $x = t^{-1}$.
By equivalent we mean that we do not obtain any new information
about the curve by considering also tropisms with first coordinate
of the opposite sign.

\begin{theorem} \label{themain}
Given a proper algebraic curve defined by a polynomial system
$\boldf(\x) = \zero$ of $n$ equations in $n$ unknowns.
Assume that the curve does not lie in a coordinate hyperplane
perpendicular to the first coordinate axis.
Then the degree of the curve corresponds to the number of
certificates times the degree of each certificate.
\end{theorem}

\noindent {\em Proof.}  If we consider the two asymptotic
homotopies in sequence, then we may lose solutions in the first move,
by putting the cutting hyperplane perpendicular to the first
coordinate axis, while in the tropisms found in the second stage
lead then to higher degrees.  The claim of the theorem is to show
that a potential loss of witness points is made up by the increase
in the degree of the Puiseux series at the very end.

Executing the asymptotic homotopies in sequence is equivalent
to combine them into one homotopy:
\begin{equation}
   \left\{
      \begin{array}{cl}
         \boldf(\x) = \zero \\
         x_1 + t ( c_0 + c_2 x_2 + \cdots + c_n x_n ) = 0 
         & \quad {\rm ~for~} t {\rm ~from~} 1 {\rm ~to~} 0.
      \end{array}
   \right.
\end{equation}
Observe that the Newton polytopes of $\boldf$ are not deformed
by this homotopy.  So the structure of the space at infinity
remains invariant as well.
Applying polyhedral endgames~\cite{HV98} (see also~\cite{MSW92}), 
the witness points will end
at roots of initial form systems.  So every witness point
corresponds to one initial root and we cannot have more
witness points than we have initial roots.

Using a weighted projective space~\cite{Cox03}, \cite{CLO98}, \cite{Ver00},
determined by the tropisms,
we can extend the initial roots to roots for $t > 0$.
So the initial roots give rise to witness points.~\qed

We end this section with a note on the complexity of the certificates.
In the best case, the initial form systems are binomial,
while in the worst case the number of monomials in the initial
form systems is of the same order of magnitude as the original system.
But even in this worst case, the number of variables drops by one
and this drop may be enough to get a more tractable problem.

\section{Tropisms and Mixed Volumes}

Degenerating witness sets is an effective but not an efficient
way to compute tropisms.  
Relating to mixed volumes, we indicate how to
compute the tropisms directly from the Newton polytopes.
Before we explain the lifting algorithm to
compute mixed volumes we point out that Gfan~\cite{Jen08}
using the algorithms of~\cite{BJSST07} is
more appropriate to compute tropisms.  Our point of relating
tropisms to mixed volumes is to investigate the connection between
the generic number of isolated roots and
the degrees of the solution curves of sparse polynomial systems.

Let $A = (A_1,A_2,\ldots,A_n)$ be the supports 
of $\boldf(\x) = \zero$.
Following~\cite{EC95}, \cite{HS95}, we summarize the mixed volume
computation in three stages:
\begin{enumerate}
\item Lift $\bfa \in A$ using a function $\omega$: 
      $\omega(A_i) \subseteq \rr^{n+1}$.
\item The facets on the lower hull of the Minkowski sum
      ${\displaystyle \sum_{i=1}^n \omega(A_i)}$ 
      spanned by one edge of each of~$\omega(A_i)$ 
      define {\em mixed cells}~$C$.
\item The mixed volume is
      ${\displaystyle V_n(A) = \sum_{\begin{array}{c}
         C \subseteq A \\ {\small C~is~mixed} \end{array}} {\rm Vol}(C)}$.
\end{enumerate}
By duality~\cite{VGC96}, 
mixed cells are defined by inner normals perpendicular
to edges of the polytopes.  
These inner normals are tropisms with positive last coordinate.

Polyhedral homotopies~\cite{HS95} 
follow the computation of the mixed cells
and the lifting function defines the powers of the new parameter~$t$.
In particular, polynomials are lifted as
\begin{equation}
  f(\x) = \sum_{\bfa \in A} c_\bfa \x^\bfa
  \quad \rightarrow \quad
  \widehat{f}(\x,t) = \sum_{\bfa \in A} c_\bfa \x^\bfa t^\omega(\bfa) 
\end{equation}
using the same lifting function~$\omega$ as before.

Now we look for solution curves of the form~(\ref{eqformsol}).
We observe that since $x_1 = t$, 
as lifting function we define $\tau$ as
$\tau(\bfa) = \deg(\x^\bfa,x_1) = a_1$ as lifting function.
The lifting~$\tau$ will work for
systems of $n-1$ equations in $n$ variables.
However, for systems with $n$ equations, 
if we use $t = x_1$, then we will have too few variables.
Therefore,
we use a slack variable $z$ in the lifting, and define $\tau$ as
\begin{equation}
   \begin{array}{ccccc}
      \tau & : & \cc[\x] & \rightarrow & \cc[z,\x] \\
           &   & \x^\bfa = x_1^{a_1} x_2^{a_2} \cdots x_n^{a_n}
           & \mapsto &
        \tau(\x^\bfa) = z^r x_2^{a_2} \cdots x_n^{a_n} t^{a_1}
   \end{array}
\end{equation}
where $r$ is some random exponent.
For exponent vectors $\bfa$, we define $\tau$ as
$\tau(a_1,a_2,\ldots,a_n) = (r,a_2,\ldots,a_n,a_1)$.
Note that $r = \omega(\bfa)$, the usual random lifting used
to calculate mixed volumes.  If we take $t = x_1$, $z = t$,
$r = \omega(\bfa)$, then we get 
$t^\omega(\bfa) x_2^{a_2} \cdots x_n^{a_n} x_1^{a_1} 
= \x^\bfa t^\omega(\bfa)$, which is the random lifting
commonly used in polyhedral homotopies~\cite{HS95}.
With the slack variable we have again as many variables
as equations and we can apply our mixed volume calculators.

We claim that the tropisms to proper algebraic curves are
in one-to-one correspondence with those inner normals to the
mixed cells of the mixed subdivision induced by the lifting~$\omega$
for which the $z$-component is zero.

\begin{proposition} \label{proptropmvc}
Let $S$ be the set of mixed cells for an $n$-tuple~$A$ in $n$-space
lifted with $\tau$, introducing a slack variable~$z$.
Then $\bfv = (\bfv_z,\bfv_\x)$ is a tropism if and only if $\bfv_z = 0$.
\end{proposition}
\noindent {\em Proof.}  Denoting the inner normals to the cells in~$S$
as $\bfv = (\bfv_z, \bfv_\x)$, points $(\bfa,\bfb)$ in each cell satisfy
\begin{equation} \label{eqmvctrop}
   r_\bfa \bfv_z  + \bfa_\x \bfv_\x = r_\bfb \bfv_z + \bfb_\x \bfv_\x,
   \quad \bfa = (r_\bfa,\bfa_\x),
   \bfb = (r_\bfb,\bfb_\x).
\end{equation}
We have to show two things:
\begin{enumerate} 
\item Each inner normal~$\bfv$ with $\bfv_z = 0$ is a tropism.

      If $\bfv_z = 0$, then 
      $\bfa_\x \bfv_\x = \bfb_\x \bfv_\x$ holds
      and the tropism is $\bfv_\x$ because $(\bfa_\x,\bfb_\x)$ span
      an edge of one of the Newton polytopes.

\item For every tropism we must have that $\bfv_z = 0$.

      Given a tropism~$\bfv_\x$, we have
      $\bfa_\x \bfv_\x = \bfb_\x \bfv_\x$.
      Rewriting~(\ref{eqmvctrop}) to solve for $\bfv_z$ leads to
\begin{equation}
   \bfv_z ( r_\bfa - r_\bfb ) = \bfb_\x \bfv_\x - \bfa_\x \bfv_\x = 0.
\end{equation}
Since $r_\bfa$ and $r_\bfb$ are random numbers, $\bfv_z = 0$.~\qed
     
\end{enumerate}

Although this is already a more efficient method than 
applying polyhedral endgames in the asymptotic homotopies
on witness sets, for practical purposes one would include
the constraint $v_z = 0$ already immediately in all feasibility tests
the mixed volume calculator does.
Ultimately, unlike the lift-and-prune approach~\cite{EC95}
for mixed volume computation, the complexity of the problem 
of computing tropisms is governed by the shape of the polytopes
and the relative position of the polytopes with respect to each other.
Tropisms lie in the common refinement of cones of inner normals
to faces of the Newton polytopes and algorithms of~\cite{BJSST07}
as implemented in Gfan~\cite{Jen08} are recommended.

Proposition~\ref{proptropmvc} relates
the number of isolated solutions to components of solutions.
Consider for example the following system:
\begin{equation}
  \boldf(\x) =
  \left\{
     \begin{array}{rcl}
        x_1 x_3 - x_1 - x_3 + 1 = 0 \\
        x_2 - 1 = 0 \\
        x_3 - \gamma = 0 & \quad & \gamma \in \cc^*.
     \end{array}
  \right.
\end{equation}
The mixed volume of $\boldf(\x) = 0$ equals one.
For all nonzero choices of~$\gamma$ (except for $\gamma = 1$),
the system will have the isolated solution~$(1,1,\gamma)$.
For $\gamma = 1$, we have the tropism $\bfv = (1,0,0)$
and the general solution line~$(x_1 = t, x_2 = 1, x_3 = 1)$.
A witness set representation for this system would need to store
two solutions: one solution for the witness set of the line
in case $\gamma = 1$ and one solution for the path leading
to the isolated root in case $\gamma \not= 1$.

\section{Preliminary Computational Experiments}

The cyclic $n$-roots problem is a widely known benchmark 
for polynomial system solvers, see e.g.~\cite{DKK03}, 
\cite{Fau99} and~\cite{LL08}.
For those~$n$ which are divisible
by a square (e.g.: $n = 4$, 8,9,12), the system is known to have
positive dimensional solution sets~\cite{Bac89}.
The general recipe to formulate
the polynomial equations for any~$n$ in this family is obvious from
the cyclic 4-roots system, given above in~(\ref{eqcyclic4}).
The permutation symmetry in cyclic $n$-roots is generated by two elements
$(x_1,x_2,x_3$, $\ldots,x_n)$ $\rightarrow$ $(x_2,x_3$, $\ldots, x_n,x_1)$
and $(x_1,x_2$, $\ldots,x_{n-1},x_n)$ 
$\rightarrow$ $(x_n,x_1,x_2,$ $\ldots, x_{n-1})$.

The numerical computations reported below were done with 
PHCpack~\cite{Ver99}, using a modified lifting in MixedVol~\cite{GLW05}. 
An alternative mixed volume calculator is DEMiCs~\cite{MT08}.
For the symbolic manipulations, Maple~11 was used.
We emphasize that the computations reported below are preliminary,
mainly to illustrate the concepts.

\subsection{cyclic 8-roots}

The program computed 29 tropisms, listed below
in Table~\ref{tab1cyc8} and Table~\ref{tab2cyc8}.

\begin{table}[hbt]
\begin{center}
$\begin{array}{r|rrr|rr|r|rrr|rrr}
 3 &  1 &  1 &  1 &  1 &  1 &  1 &  1 &  1 &  1 &  1 &  1 &  1 \\
-1 &  1 &  1 & -3 &  1 & -1 & -1 & -1 &  0 & -1 &  0 & -1 & -1 \\
-1 &  1 & -3 &  1 & -1 & -1 &  1 &  1 & -1 &  1 & -1 &  1 &  0 \\
-1 & -3 &  1 &  1 & -1 &  1 & -1 &  0 &  0 & -1 &  1 & -1 &  1 \\
 3 &  1 &  1 &  1 &  1 &  1 &  1 & -1 &  1 &  1 & -1 &  0 &  0 \\
-1 &  1 &  1 & -3 &  1 & -1 & -1 &  0 & -1 &  0 &  1 &  1 & -1 \\
-1 &  1 & -3 &  1 & -1 & -1 &  1 &  1 &  1 & -1 & -1 &  0 &  1 \\
-1 & -3 &  1 &  1 & -1 &  1 & -1 & -1 & -1 &  0 &  0 & -1 & -1 \\
\end{array}$
\caption{First list of 13 tropisms for cyclic 8-roots, separated in 6 orbits.}
\label{tab1cyc8}
\end{center}
\end{table}

\begin{table}[hbt]
\begin{center}
$\begin{array}{rr|rr|rr|rr|rr|rr|rr|rr}
 1 &  1 &  1 &  1 &  1 &  1 &  1 &  1 &  1 &  1 &  1 &  1 &  1 &  1 &  1 &  1 \\
 0 & -1 &  0 & -1 &  0 & -1 &  0 & -1 &  0 & -1 &  0 & -1 &  0 & -1 &  0 &  0 \\
 0 &  1 &  0 &  0 &  0 &  1 &  0 &  0 &  0 &  0 & -1 &  1 & -1 &  0 & -1 & -1 \\
 0 &  0 &  0 &  1 & -1 &  0 & -1 &  1 & -1 &  0 &  0 &  0 &  1 &  0 &  0 &  1 \\
-1 &  0 & -1 &  0 &  0 &  0 &  0 &  0 &  1 &  1 &  0 & -1 & -1 &  0 &  0 &  0 \\
 0 &  0 &  1 &  0 &  0 & -1 &  1 &  0 & -1 &  0 &  0 &  0 &  0 &  1 &  1 & -1 \\
 1 & -1 & -1 &  0 &  1 &  0 & -1 & -1 &  0 &  0 &  1 &  0 &  0 &  0 &  0 &  0 \\
-1 &  0 &  0 & -1 & -1 &  0 &  0 &  0 &  0 & -1 & -1 &  0 &  0 & -1 & -1 &  0
\end{array}$
\end{center}
\caption{Second list of 16 tropisms for cyclic 8-roots, separated in 8 orbits.}
\label{tab2cyc8}
\end{table}

The last tropism from Table~\ref{tab2cyc8}:
$\bfv = (1,0,-1,1,0,-1,0,0)$, the initial form system is
\begin{equation}
  {\rm in}_\bfv \boldf(\x) =
  \left\{
     \begin{array}{r}
        x_3 + x_6 = 0 \\
        x_2 x_3 + x_5 x_6 + x_6 x_7 = 0 \\
        x_5 x_6 x_7 + x_6 x_7 x_8 = 0 \\
        x_3 x_4 x_5 x_6 + x_5 x_6 x_7 x_8 = 0 \\
        x_2 x_3 x_4 x_5 x_6 + x_3 x_4 x_5 x_6 x_7  = 0 \\
        x_2 x_3 x_4 x_5 x_6 x_7 + x_3 x_4 x_5 x_6 x_7 x_8
        + x_6 x_7 x_8 x_1 x_2 x_3 = 0 \\
        x_2 x_3 x_4 x_5 x_6 x_7 x_8 + x_5 x_6 x_7 x_8 x_1 x_2 x_3 = 0 \\
        x_1 x_2 x_3 x_4 x_5 x_6 x_7 x_8 - 1 = 0
     \end{array}
  \right.
\end{equation}
To reduce the system into a simpler form, we perform the coordinate
transformation, dividing out $y_1$ in the initial form system:
\begin{equation}
  \left\{
     \begin{array}{l}
        x_1 = y_1 \\
        x_2 = y_1^0 y_2 \\
        x_3 = y_1^{-1} y_3 \\
        x_4 = y_1^1 y_4 \\
        x_5 = y_1^0 y_5 \\
        x_6 = y_1^{-1} y_6 \\
        x_7 = y_1^0 y_7 \\
        x_8 = y_1^0 y_8 \\
     \end{array}
  \right.
  \quad
  {\rm in}_\bfv \boldf(\y) = 
  \left\{
     \begin{array}{r}
        y_3 + y_6 = 0 \\
        y_2 y_3 + y_5 y_6 + y_6 y_7 = 0 \\
        y_5 y_6 y_7 + y_6 y_7 y_8 = 0 \\
        y_3 y_4 y_5 y_6 + y_5 y_6 y_7 y_8 = 0 \\
        y_2 y_3 y_4 y_5 y_6 + y_3 y_4 y_5 y_6 y_7 = 0 \\
        y_2 y_3 y_4 y_5 y_6 y_7 + y_3 y_4 y_5 y_6 y_7 y_8
        + y_6 y_7 y_8 y_2 y_3 = 0 \\
        y_2 y_3 y_4 y_5 y_6 y_7 y_8 + y_5 y_6 y_7 y_8 y_2 y_3 = 0 \\
        y_2 y_3 y_4 y_5 y_6 y_7 y_8 - 1 = 0 
     \end{array}
  \right.
\end{equation}
To solve the overconstrained initial form, we introduce
a slack variable~$s$ and generate eight random numbers~$\gamma_k \in \cc^*$,
$k=1,2,\ldots,8$ to multiply~$s$ with.  Then we solve
\begin{equation}
  {\rm in}_\bfv \boldf(\y,s) = 
  \left\{
     \begin{array}{r}
        y_3 + y_6 + \gamma_1 s = 0 \\
        y_2 y_3 + y_5 y_6 + y_6 y_7 + \gamma_2 s = 0 \\
        y_5 y_6 y_7 + y_6 y_7 y_8 + \gamma_3 s = 0 \\
        y_3 y_4 y_5 y_6 + y_5 y_6 y_7 y_8 + \gamma_4 s = 0 \\
        y_2 y_3 y_4 y_5 y_6 + y_3 y_4 y_5 y_6 y_7 + \gamma_5 s = 0 \\
        y_2 y_3 y_4 y_5 y_6 y_7 + y_3 y_4 y_5 y_6 y_7 y_8
        + y_6 y_7 y_8 y_2 y_3 + \gamma_6 s = 0 \\
        y_2 y_3 y_4 y_5 y_6 y_7 y_8 + y_5 y_6 y_7 y_8 y_2 y_3
        + \gamma_7 s = 0 \\
        y_2 y_3 y_4 y_5 y_6 y_7 y_8 - 1 + \gamma_8 s = 0 
     \end{array}
  \right.
\end{equation}
The mixed volume for this system equals 25 and is exact.
Of the 25 solutions, eight solutions have $s=0$ and are thus proper
solutions.  Two of the eight solutions are real.
Denoting $I = \sqrt{-1}$, we select the solution
\begin{equation}
   y_2 = -\frac{1}{2} + \frac{I}{2}, y_3 = -I, y_4 = -1,
   y_5 = -1 + I, y_6 = I, y_7 = \frac{1}{2} - \frac{I}{2},
   y_8 = 1 - I.
\end{equation}
Then we look at the first term of the series
\begin{equation}
   \left\{
      \begin{array}{lllcl}
         y_1 & = & t \\
         y_2 & = & -\frac{1}{2} + \frac{I}{2} & + & z_2~t\\
         y_3 & = & -I & + & z_3~t\\
         y_4 & = & -1 & + & z_4~t\\
         y_5 & = & -1 + I & + & z_5~t\\
         y_6 & = & I & + & z_6~t\\
         y_7 & = & \frac{1}{2} - \frac{I}{2} & + & z_7~t\\
         y_8 & = & 1 - I & + & z_8~t.
      \end{array}
   \right.
\end{equation}
To decide whether the solution is isolated or not, we need
to find values for the coefficient of the second term in
the expansion.  Substituting the series in the system 
in $\boldf(\y) = \zero$ and selecting the lowest order terms
in~$t$ leads to an overdetermined {\em linear} system
in the $z_k$ variables.  Solving with Maple yields
\begin{equation}
   z_2 = -\frac{1}{2}, z_3 = -1 + I, z_4 = 0, z_5 = -1,
   z_6 = 1-I, z_7 = \frac{1}{2}, z_8 = 1.
\end{equation}
Substituting the series in $\boldf(\y)$ and we see that the
result is~$O(t^2)$.

\subsection{cyclic 12-roots}

Extrapolating on the tropism for cyclic 4-roots, 
we considered $\bfv = (+1,-1,+1,-1$, $+1,-1,+1,-1$, $+1,-1,+1,-1)$.
For this tropism
the first term of the Puiseux series expansion

\begin{equation}
\begin{array}{lcl}
   x_1 = t & \quad
 & x_2 = t^{-1} \left(  \frac{1}{2} - \frac{1}{2} i \sqrt{3} \right) \\
\vspace{-3mm} \\
   x_3 = -t &
 & x_4 = t^{-1} \left( -\frac{1}{2} - \frac{1}{2} i \sqrt{3} \right) \\
\vspace{-3mm} \\
   x_5 = t      \left( -\frac{1}{2} + \frac{1}{2} i \sqrt{3} \right)  &
 & x_6 = t^{-1} \left(  \frac{1}{2} + \frac{1}{2} i \sqrt{3} \right) \\
\vspace{-3mm} \\
   x_7 = - t &
 & x_8 = t^{-1} \left( -\frac{1}{2} + \frac{1}{2} i \sqrt{3} \right) \\
\vspace{-3mm} \\
   x_9 = t &
 & x_{10} = t^{-1} \left(  \frac{1}{2} + \frac{1}{2} i \sqrt{3} \right) \\
\vspace{-3mm} \\
   x_{11} = t      \left(  \frac{1}{2} - \frac{1}{2} i \sqrt{3} \right) &
 & x_{12} = t^{-1} \left( -\frac{1}{2} - \frac{1}{2} i \sqrt{3} \right) 
\end{array}
\end{equation}
makes the system entirely and exactly equal to zero.
Because of the symmetry, we have
five other solution curves of this type.
This is an exact certificate that shows cyclic 12-roots has
a curve of degree two.   Note that 0.866025403784439 is
close enough to $\sqrt{3}/2$ for us to recognize.
The numerical determination of algebraic numbers in general
is done via the integer relation detection algorithm of~\cite{BF89},
see also~\cite{Bai00}.

The initial root that led to this certificate was one of 
roots of an initial form system with mixed volume --- after
adding one slack variable to make the system square ---
equal to~49,816.  Note that this number is much less than
the mixed volume of the original system:~500,352.

\section{Conclusions and Future Directions}

In this paper concepts of numerical algebraic geometry were applied
--- witness sets and endgames ---
to determine the orientation of tropisms which could lead
to certificates for proper algebraic curves.

We sketched how a polyhedral method could pick up
all proper algebraic curves defined by a polynomial system.
To extend this to solution sets of any dimension, say~$k$, one
would need to consider Puiseux series with $k$ free parameters
and look for $k$ dimensional cones of tropisms.
This generalization leads the development of multiparameter
polyhedral homotopies which may be of independent interest
to numerical analysis.

\bibliographystyle{plain}

\end{document}